\documentclass[11pt]{amsart}
\usepackage{amssymb,pstricks,pst-plot}

\input xy
\xyoption{all}

\newcommand{\C}{\mathbb{C}}
\newcommand{\D}{\mathbb{D}}
\newcommand{\T}{\mathbb{T}}
\newcommand{\clD}{{\overline{\,\D}}}

\newcommand{\Z}{\mathbb{Z}}

\newcommand{\R}{\mathbb{R}}

\newcommand{\bT}{\mathbb{T}}

\newcommand\E{\mathrm{e}}
\newcommand\I{\mathrm{i}}

\newcommand{\cA}{\mathcal{A}}

\newcommand{\cC}{\mathcal{C}}

\newcommand{\cL}{\mathcal{L}}

\newcommand{\ol}{\overline}

\newcommand\bs{\backslash}

\newcommand\di{\partial}

\newcommand\dibar{\bar\partial}

\newtheorem{theorem}{Theorem}[section]
\newtheorem{lemma}[theorem]{Lemma}

\theoremstyle{definition}

\numberwithin{equation}{section}

\def\ss{\Subset}

\begin{document}
\title{A disc formula for plurisubharmonic subextensions in manifolds}        
\author{Barbara Drinovec Drnov\v sek}
\address{Faculty of Mathematics and Physics, University of Ljubljana, 
and Institute of Mathematics, Physics and Mechanics, Jadranska 19, 
1000 Ljubljana, Slovenia}
\email{barbara.drinovec@fmf.uni-lj.si}
\thanks{Research supported by grant P1-0291, Republic of Slovenia.}

%
%
%
%

\date{\today} 

%
%
%
%
\begin{abstract}
  We provide a sufficient condition for open sets $W$ and $X$ such that
  a disc formula for the largest plurisubharmonic subextension of an upper       
  semicontinuous function on a domain $W$ to a complex  manifold  $X$ holds.
\end{abstract}

\maketitle

%
%
%
%
%
%

\section{Introduction and the main result}
\label{intro}
The theory of disc functionals was initiated by 
E. A. Poletsky in the late 1980s \cite{Poletsky1991} and it
offers a different approach to certain 
extremal functions of pluripotential theory
(see Klimek \cite{Klimek}).
Extremal functions are usually 
defined as suprema of classes of plurisubharmonic functions 
with certain properties, and many of them are 
envelopes of appropriate disc functionals: 
the largest plurisubharmonic minorant of an 
upper semicontinuous function,
the largest nonpositive plurisubharmonic function 
whose Levi form is bounded below by the 
Levi form of a given plurisubharmonic function,
pluricomplex Green functions and
Siciak-Zaharyuta extremal function  \cite{Poletsky1993,Bu-Schachermayer,Larusson-Sigurdsson1998,
 Rosay1,Rosay2,Larusson-Sigurdsson2003, Rashkovskii-Sigurdsson2005,
 Larusson-Sigurdsson2005,Larusson-Sigurdsson2007, Magnusson-Sigurdsson2007,Larusson-Sigurdsson2009,
 Magnusson2011,DF-Indiana, DF-Polonici, Kuzman}.

Recently, F. L\'arusson and E. A. Poletsky \cite{Larusson-Poletsky2012} proved a disc formula for the largest 
plurisubharmonic subextension of an upper semicontinuous function on a domain $W$  to a larger domain $X$ in a Stein manifold under suitable conditions on $W$ and $X$.
Their conditions depend on the topological properties of the space $\cA_X^W$ of analytic discs in $X$ with boundaries in $W$. 
They also provide an example which shows that the disc formula can fail in general.
The aim of this note is to provide a more geometric condition on $X$ and $W$ which 
also implies the disc formula. 

\begin{theorem}
\label{main}
Let $X$ be a $(n-1)$-complete complex manifold of dimension $n>1$,
with smooth $(n-1)$-complete exhaustion function $\rho\colon X\to \R$.
Let $W=\{x\in X\colon\rho(x)> c\}$ be a superlevel set of $\rho$ for some $c\in\R$.
If $\phi\colon W\to [-\infty, \infty)$  is upper semicontinuous, then
$$\sup\{u\in PSH(X)\colon u|W\le \phi\}=\inf\left\{\int_\T \phi\circ f\,d\lambda\colon f\in \cA_X^W,\ f(0)=x\right\}.$$ 
\end{theorem}

First we set some notations. Let $\D=\{\zeta\in \C\colon |\zeta|< 1\}$ and $\bT=b\D=\{\zeta\in\C\colon |\zeta|=1\}$. 
Given a domain $W$ in a complex manifold $X$ we denote by $\cA_X$ 
the set of all {\em analytic discs} in $X$, that is, continuous maps $\clD \to X$ that are holomorphic in $\D$, and by $\cA_X^W$ the set of all
analytic discs in $X$ with boundaries in $W$. We denote by $PSH(X)$ the class of
all plurisubharmonic functions on $X$ (we take the constant $-\infty$ to be plurisubharmonic).
If $\phi\colon W\to [-\infty, \infty)$ is an upper semicontinuous function, then 
 we denote by $S\phi$ the supremum of all 
plurisubharmonic subextensions of $u$, that is plurisubharmonic functions $u\in PSH(X)$
satisfying $u|W\le \phi$. 
As noted in \cite{Larusson-Poletsky2012}
$S\phi$ is plurisubharmonic if, for example, $X$ is covered by analytic discs 
with boundaries in $W$. In this case, $S\phi$ is 
the largest plurisubharmonic subextension of $\phi$.
The {\em Poisson functional} associated to $\phi$
is the map $H_\phi\colon \cA_X^W \to[-\infty, \infty)$  defined by
\[
	H_\phi(f)= \int_\T \phi\circ f\,d\lambda,\qquad f\in \cA_X^W,
\]
where $\lambda$ denotes the normalized arc length measure on $\T$.
As in \cite{Larusson-Poletsky2012} we define
the {\em Poisson envelope} of $\phi$ with respect to $\cA_X^W$ by 
$$E_{\cA_X^W}\phi(x)=\inf\{H_\phi(f)\colon f\in \cA_X^W,\ f(0)=x\}.$$

Recall that a smooth function 
$\rho\colon X\to\R$ on a complex manifold $X$ is said to be {\em $q$-convex} 
on an open subset $U\subset X$ (in the sense of Andreotti-Grauert \cite{AG},
\cite[Definition\ 1.4, p.\ 263]{Grauert2})
if its Levi form $i\di\dibar \rho=\cL_\rho$ 
has at most $q-1$ negative or zero eigenvalues at each point of $U$.
The manifold $X$ is {\em $q$-complete}, if it admits 
a smooth exhaustion function $\rho\colon X\to\R$ which is $q$-convex on $X$.
A 1-complete complex manifold is just a Stein manifold.

Note that if $X$ and $W$ satisfy the conditions in Theorem \ref{main}, 
then $X$ is covered by analytic discs with
boundaries in $W$ \cite{BD-Fourier,DF-Duke}.

%
%
%
%
%
%
%

\section{Proof}

In the proof we use the method of gluing sprays, which was developed in 
\cite{DF-Duke}, and it was used in the Poletsky theory in 
\cite{DF-Indiana,DF-Polonici} (For an exposition of the method 
we refer to \cite[\S 5.8--\S 5.9]{F-book}; 
a brief exposition can also be found in \cite{DF-Indiana}).

The following lemma is essentially \cite[Lemma 6.2]{DF-Duke}: 

\begin{lemma}
\label{smallstep}
Let $X$ be a complex manifold of  dimension $n>1$, and let $d$ be a complete 
distance function on $X$.
Let $A\Subset  X$ be relatively compact open subset of $X$ and 
let $B$ be a 2-convex bump on $A$ \cite[Definition 6.1]{DF-Duke}.
Let $P$ be a domain in $\C^N$ containing $0$ and 
assume that $f\colon P\times\clD\to X$ is a spray of discs
with the exceptional set containing $0$ 
such that $f_0( b\D)\cap \bar A=\emptyset$. 
(Here $f_0=f(0,\cdotp)$ is the core map of the spray.)
Let $E \subset \T$ be a union of finitely many closed
circular arcs  such that 
$f_0(\overline{\T\setminus E})\cap  \overline{A\cup B}=\emptyset$.

Given $\epsilon>0$  there are a domain $P'\subset P$ containing $0\in\C^N$ and 
a spray  of maps $g\colon P'\times\clD\to X$ 
 with the exceptional set containing $0$
such that 
\begin{itemize}
\item[(i)]   $g(t,z)\cap \overline{A\cup B}=\emptyset$ for $t\in P'$ and $z\in\T$,
\item[(ii)] $d(g(t,z),f(t,z))<\epsilon$ for $t\in P'$ and for 
$z\in \overline{\T\setminus E}$, and
\item[(iii)] $f_0(0)=g_0(0)$.
\end{itemize}
\end{lemma}

Let us give the idea of the proof of the above lemma.
A 2-convex bump
gives a continuously varying family of small analytic discs in $X$ with boundaries 
outside the bump along the part of the boundary of the initial disc that is mapped
into the bump.
By solving a certain Riemann-Hilbert boundary value problem we obtain a 
continuous map defined on a neighborhood of a boundary arc in the closed unit disc,
holomorphic in the interior, such that the boundary arc 
is mapped close to the boundary of the above chosen family of analytic discs, and hence the boundary arc is  mapped outside the bump. 
The globalization of the construction is obtained by the method of patching 
holomorphic sprays. The local modification method gives a new spray 
near a part of the boundary $\T$; by insuring that the two sprays are sufficiently 
close to each other on the intersection of their domains,  
we patch them into a new spray. 
We obtain property (ii) by possibly shrinking the parameter set $P$ 
at the start. 

We use sprays that are continuous up to the
boundary and as it is remarked in \cite{DF-Duke}, 
the assumption that sprays are of class $\cC^r$, $r\ge 2$,
is needed only for the existence of a Stein neighborhood. 
Using \cite[Corollary  1.3]{FFAsian} instead of \cite[Theorem 2.6]{DF-Duke} we obtain
the same result for all $r\in\Z_+$. Here we shall use these results with $r=0$.

The following lemma is essentially obtained by the same proof as
\cite[Lemma 6.3]{DF-Duke}: property (ii) below 
is a consequence of Lemma \ref{smallstep} (ii) which we use in the proof
instead of using \cite[Lemma 6.2]{DF-Duke}. 

\begin{lemma}
\label{bigstep}
Let $X$ be a $(n-1)$-complete complex manifold of dimension $n>1$,
with smooth $(n-1)$-complete exhaustion function $\rho\colon X\to \R$ 
and let $d$ be a complete distance function on $X$.
Let $f\in \cA_X$ be an analytic disc, let $c_0<c_1$ and assume that 
$\rho(f(z))<c_1$ for all $z\in \T$ and furthermore assume that $E \subset \T$ is
a union of finitely many closed circular arcs such that 
$\rho(f(z))\in(c_0,c_1)$ for all $z\in\overline{\T\setminus E}$.
Then for each $\eta>0$ there exists 
an analytic disc $g\in \cA_X$  
satisfying the following properties:
\begin{itemize}
\item[(i)]   $\rho(g(z))\in(c_0,c_1)$ for $ z\in \T$,
\item[(ii)]  $d(g(z),f(z))<\eta$ for $z\in \overline{\T\setminus E}$, and
\item[(iii)]  $f(0)=g(0)$.
\end{itemize}
\end{lemma}

Let us review the idea of the proof.
Since the Levi form $\cL_\rho$ is assumed to have at least two positive eigenvalues 
at every point $x_0\in X$, we get at least one positive eigenvalue in a direction 
tangential to the regular level set of $\rho$ at each point $x_0$. So we can 
choose local coordinates near $x_0$ such that the sublevel set 
$\{x\colon \rho(x)\le \rho(x_0)\}$ near $x_0$
is convex in one direction and  we can fill the space between 
regular sublevel sets by adding finitely many 2-convex bumps.

By Lemma \ref{smallstep} we can push the boundary of the analytic disc outside 
a $2$-convex bump while we make only a small perturbation on the part of the boundary
that already lies outside some neighborhood of the bump. In finitely many steps
we can push the boundary to a given higher superlevel set of $\rho$. 
To cross each critical level set of $\rho$ we use a 
different function constructed especially for this purpose: we make a small 
perturbation if necessary and we proceed by the noncritical 
procedure with perturbed $(n-1)$-convex function. Once the boundary is above the 
critical level set we proceed with the noncritical construction for 
the original exhaustion function.

\textit{Proof of Theorem \ref{main}.}
Note that $E_{\cA_X^W}\phi\le \phi$ on $W$ since constant discs in $W$ belong
to $\cA_X^W$. Once we prove that $E_{\cA_X^W}\phi$ is plurisubharmonic the equality 
in the theorem holds.

The proof is reduced to the case when the function 
$\phi\colon W\to\R$ is continuous and bounded from below 
exactly as in the case of Poisson envelope \cite{DF-Indiana}.


Next we show that $E_{\cA_X^W}\phi$ is upper semicontinuous
on $X$. The proof goes along the same lines as in \cite{DF-Indiana};
we include it for the convenience of the reader since it applies
sprays.
Pick a point $x\in X$ and a number $\epsilon>0$. 
Assume first that $E_{\cA_X^W}\phi(x)>-\infty$.
By the definition of $E_{\cA_X^W}\phi$ there exists a disc 
$f_0\in \cA_X^W$ with $f_0(0)=x$ such that 
$E_{\cA_X^W}\phi(x)\le H_\phi(f_0) < E_{\cA_X^W}\phi(x)+\epsilon$.
We embed $f_0$ as the central map $f_0=f(0,\cdotp)$ in a 
spray of holomorphic discs $f\colon  P\times\clD\to X$,
where $P$ is an open set in $\C^m$ containing the origin.
If $P'$, $0\in P'\subset P$, is small enough then
$f(t,\cdotp)\in \cA_X^W$, $H_\phi(f(t,\cdotp))< H_\phi(f_0)+\epsilon$ 
for each $t\in P'$, and hence 
\[
	E_{\cA_X^W}\phi(f(t,0))\le H_\phi(f(t,\cdotp)) \le  H_\phi(f_0)+\epsilon < E_{\cA_X^W}\phi(x)+2\epsilon.
\]
By the domination property the set 
$\{f(t,0)\colon t\in P'\}$ fills a neighborhood of the point 
$x=f_0(0)$ in $X$, so we see that $E_{\cA_X^W}\phi$ is upper semicontinuous at $x$.
A similar argument works at points where $E_{\cA_X^W}\phi(x)=-\infty$.

The main part is to prove that $E_{\cA_X^W}\phi$ is plurisubharmonic on $X$:
we need to show that for every analytic disc 
$h\in \cA_X$  we have the submeanvalue property
\begin{equation*}
 E_{\cA_X^W}\phi(h(0)) \le \int_\T (E_{\cA_X^W}\phi)(h(w))d\lambda(w)\,.
\end{equation*}
Therefore, it is enough to construct for each $\epsilon>0$ and for each
continuous function $v\ge E_{\cA_X^W}\phi$ on $X$, an analytic disc
$g\in \cA_X^W$ such that $g(0)=h(0)$, and 
\begin{equation}
\label{eq:submean}
H_\phi(g)\le H_v(h)+\epsilon .
\end{equation}

We proceed in two steps. First we construct an analytic disc $f$ with 
large part of the boundary lying in $W$ and satisfying an 
integral estimate similar to (\ref{eq:submean}) on the part of the boundary
lying in $W$. Then we correct the disc around small pieces of the boundary
while we do not perturb the rest of the boundary too much and we obtain $g$.

More precisely, for any given $\epsilon>0$
we shall construct $f\in\cA_X$, a union of finitely many closed
circular arcs $E \subset \T$, and real numbers $c<c_0<c_1$ 
with the following properties:
\begin{itemize}
\item [(a)] $\rho(f(w))<c_1$ for every $w\in\T$,
\item [(b)] $\rho(f(w))\in (c_0,c_1)$ for every $w\in\overline{\T\setminus E}$,
\item[(c)] $\int_{\T\setminus E}\phi\circ f\,d\lambda\le H_v(h)+2\epsilon/3$,
\item[(d)] $\max\{\phi(x)\colon \rho(x)\in [c_0,c_1]\}\cdot\lambda(E)\le \epsilon/6$,
\item [(e)]$f(0)=h(0)$.
\end{itemize}

Fix a point $w \in \bT$.
By the definition of $v$ there exists an analytic disc 
$g_w  \in \cA_X^W$ with center $g_w(0)=h(w)$
such that
\begin{equation}
\label{eq:suboptimal}
	\int_\T \phi\bigl(g_w(z)\bigr)\, d\lambda(z)
	 < v(h(w)) +\epsilon/6.
\end{equation}
We embed the disc $g_w$ into a dominating spray of discs 
$G(t,\cdotp)\in \cA_X$ depending holomorphically on the 
parameter $t$ in some ball $P$ in Euclidean space $\C^m$ 
with the central map $G(0,\cdotp)=g_w$. 
Since the spray $G$ is dominating, the set 
$G(P,0)=\{G(t,0)\colon t\in P\}$ 
covers a neighborhood of the point $f(w)$ in $X$.
By shrinking the parameter set $P$ if necessary we obtain, in 
addition, that $G(t,\cdotp)\in \cA_X^W$ for every $t\in P$.
By the implicit mapping theorem there is a 
disc $D\subset r_0\D$ ($r_0>1$) centered at the point $w\in \bT$ 
and a holomorphic map $\varphi\colon D \to P$ such that 
\[	
	\varphi(w)= 0 \quad {\rm and}\quad 
	G(\varphi(\zeta),0)=h(\zeta),\quad \zeta\in D.
\]
The map $g\colon D\times \clD\to X$ defined by  
\[
	g(\zeta,z)=G(\varphi(\zeta),z), \qquad \zeta\in D,\ z\in\clD
\]
is continuous, holomorphic in $D\times\D$, and 
\[
   g(w,\cdotp)=G(0,\cdotp)=g_w; \qquad 
   g(\zeta,0)= h(\zeta),\quad \zeta\in D.
\]
Since $g(\zeta,\cdotp)$ is uniformly close to 
$g_w$ when $\zeta$ is close to $w$,
it follows from (\ref{eq:suboptimal}) that 
there is a small arc $I\Subset \bT\cap D$ 
around the point $w$ such that
\[
      \int_{I} d\lambda(w)\!  \int_\T 
       \phi\bigl(g(w,z)\bigr) \, 
       d\lambda(z)  
       <
       \int_I v\bigl(h(w)\bigr)\,d\lambda(w)  + 
       \frac{\lambda(I)\,\epsilon}{6}.     
\]
By repeating this construction at other points of $\bT$
we find finitely many open circular arcs $I''_j$ and 
discs $D_j$ contained in $r_0\D$ such that 
$I''_j \ss \bT\cap D_j$ $(j=1,\ldots,l)$, $\cup I''_j=\T$,
and holomorphic families of discs $g_j(\zeta,z)$ for $\zeta\in \bar D_j$
and $z\in\clD$ such that 
%
\begin{equation}
\label{est1}
    \int_{I_j} d\lambda(w)\!  \int_\T 
       \phi\bigl(g_j(w,z)\bigr) \, 
       d\lambda(z)  
       < 
       \int_{I_j} v\bigl(h(w)\bigr)\,d\lambda(w)  + 
       \frac{\lambda(I_j)\,\epsilon}{6}.
\end{equation}
for any $I_j\subset I''_j$.

There are real numbers $c<c_0<c_1$ such that $\rho(x)\in(c_0,c_1)$
for each $x \in\bigcup g_j(\ol D_j,\clD)$.
Choose open circular arcs $I_j$ and $I'_j$ such that
$I_j\ss I'_j\ss I''_j$,  $\ol I'_j\cap \ol I'_k=\emptyset$ if $j\ne k$ and
the set $E=\bT\bs \cup_{j=1}^m I_j$ has so small measure that
$\lambda(E)\cdot \max\{\phi(x)\colon \rho(x)\in [c_0,c_1]\}\le \epsilon/6$
(and satisfies (d)), and 
\begin{equation}
\label{estx}
\left|\int_{E}v(h(w))d\lambda (w)\right|<\epsilon/6.
\end{equation}

For each $j=1,\ldots,l$ we choose a smoothly bounded simply 
connected domain $\Delta_j \subset D_j\cap \D$ 
such that $\ol\Delta_j \subset D_j$, 
$\ol\Delta_j\cap \bT=\ol{I'_j}$, and  
$\ol \Delta_j\cap \ol\Delta_k=\emptyset$ when  
$1\le j\ne k\le l$. 
Let $\chi\colon \C\to [0,1]$ be a smooth function
such that $\chi=1$ on $\cup_{j=1}^l I_j\subset \bT$ 
and $\chi=0$ in a neighborhood of the set 
$\clD \setminus \cup_{j=1}^l (\Delta_j\cup I'_j)$. 
Consider the map $\xi\colon\clD \times \clD\to X$ 
defined by
\begin{equation}
\label{discs-xi}
	\xi(\zeta,z)= 
		\begin{cases}  
					g_j\bigl(\zeta,\chi(\zeta)z\bigr), 
					  & \zeta \in \ol\Delta_j,\ z\in \clD,\ j=1,\ldots,l; \cr
  			h(\zeta),  & \chi(\zeta)=0,\ z\in \clD. 
  \end{cases}
\end{equation}
%
Note that $\xi$ is continuous and is holomorphic 
in the second variable. Then 
%
\begin{equation}
\label{est2}
   \int_{\T\setminus E}d\lambda(w)\!  \int_\T 
       \phi\bigl(\xi(w,z)\bigr) \, 
        d\lambda(z)
   < \int_\T  v\bigl(h(w)\bigr)\,d\lambda(w)
    + \epsilon/3.
\end{equation}
Indeed, the integral over $\cup_{j=1}^l I_j$
is estimated by adding up the inequalities (\ref{est1}) and the rest by 
(\ref{estx}).

Fix an index $j\in\{1,\ldots,l\}$.
We shall apply \cite[Lemma 3.1]{DF-Indiana}  over $\ol\Delta_j$ 
to find an analytic disc $f'_j \colon \ol \Delta_j\to X$ 
that approximates $h$ uniformly as close as desired outside of a 
small neighborhood of the arc $\ol I_j$ 
in $\ol \Delta_j$ and  satisfies  the estimate
%
\begin{equation}
\label{local1}
	\int_{I_j} 
	 \phi \bigl(f'_j (w)\bigr) \,d\lambda(w) < 
	\int_{I_j}d\lambda(w) \!\int_\T  
       \phi\bigl( \xi(w,z)\bigr) \, 
        d\lambda(z)
       + \frac{\lambda(I_j)\,\epsilon}{6}.
\end{equation}
Consider the function 
\[
	\phi'_j(\zeta,z)=\phi\bigl(g_j\bigl(\zeta,z)\bigr),\qquad
	\zeta\in \ol\Delta_j,\ z\in\clD
\]
and the smooth family of analytic discs in $\C^2_{(\zeta,z)}$ 
given by
\[
	g'_j(\zeta,z)=\bigl(\zeta,\chi(\zeta)\, z \bigr), \qquad 
	\zeta\in b\Delta_j,\ z\in\clD.
\]
Recall from (\ref{discs-xi}) that 
$\xi(\zeta,z)= g_j\bigl(\zeta,\chi(\zeta)z\bigr)$
for $\zeta \in \ol\Delta_j$ and $z\in \clD$ thus
\begin{equation}
\label{xi}
	\xi=g_j\circ g'_j\quad{\rm and}\quad \phi'_j=\phi\circ g_j
	\quad{\rm on}\quad b\Delta_j\times\clD.
\end{equation}
Applying \cite[Lemma 3.1]{DF-Indiana} with $\D$ replaced 
by $\Delta_j$, $u$ replaced by $\phi'_j$ and 
$g$ replaced by $g'_j$ furnishes an analytic disc
$h'_j\colon \ol \Delta_j \to \C^2$ which approximates 
the disc $\zeta\mapsto (\zeta,0)$ outside of a 
small neighborhood of the arc $\ol I_j$,
%
\begin{equation}
\label{local2}
	\int_{I_j} 
	\phi'_j \bigl(h'_j (w)\bigr) \,d\lambda(w) <
	\int_{I_j}d\lambda(w) \!\int_\T  
       \phi'\bigl( g'_j(w,z)\bigr) \, 
     d\lambda(z)
        +  \frac{\lambda(I_j)\,\epsilon}{6},
\end{equation}
and for every $\E^{\I t}\in I_j$ the point $h_j'(\E^{\I t})$ is 
so close to the set $\{\E^{\I t}\}\times\T$ that $\rho(g_j(h_j'(\E^{\I t})))\in(c_0,c_1)$.
Since $g_j\colon \ol\Delta_j\times\clD\to X$ 
is holomorphic in the interior $\Delta_j\times\D$,  
the map
\begin{equation}
\label{fjprime}
	f'_j:=g_j\circ h'_j\colon \ol\Delta_j\to X
\end{equation}
is an analytic disc in $X$ such that 
$f'_j(\zeta)\approx g_j(\zeta,0)=f(\zeta)$
for $\zeta$ outside a small neighborhood of $\ol I_j$
in $\ol \Delta_j$. From (\ref{xi}) and (\ref{fjprime})
it follows that
\[
	\phi\circ f'_j= \phi\circ g_j \circ h'_j = \phi'_j\circ h'_j,
	\qquad
	\phi\circ\xi = \phi\circ g_j\circ g'_j = \phi'_j\circ g'_j
\]
hold on $b\Delta_j\times \clD$.
Hence the integrals in (\ref{local1}) equal
the corresponding integrals in (\ref{local2}), and so the disc $f'_j$
satisfies the desired properties.

If the approximation of $h$ by $f'_j$ is close enough 
for each $j=1,\ldots,l$, we can glue this collection of discs into a single analytic disc 
$f\colon \ol{\,\D}\to X$ with the same center as $h$, which approximates $h$ away from 
the union of arcs $\cup_{j=1}^l I_j$, and which approximates 
the disc $f'_j$ over a neighborhood of $\ol I_j$ for each $j$
exactly as in the proof of \cite[Theorem 1.1]{DF-Indiana}. 
If the approximation is good enough the properties (a) and (b) hold.
In particular, we can insure that for each $j=1,\ldots, l$ we have
\[
	\int_{I_j} \phi\bigl(f(w)\bigr) \,d\lambda(w)<
	\int_{I_j}d\lambda(w) \!\int_\T  
       \phi\bigl(g(w,z)\bigr) \, 
   d\lambda(z)+\frac{\lambda(I_j)\,\epsilon}{3}.
\]

By adding up these terms  and  by the inequality (\ref{est2}) 
we get
%
\begin{align}
\label{small-increase3}
	\int_{\T\setminus E} \phi(f(w)) \, d\lambda(w)
	&< \int_{\T\setminus E} d\lambda(w)\! \int_\T 
       \phi\bigl(\xi(w,z)\bigr) \, 
        d\lambda(z)+ \epsilon/3 \nonumber\\
      &< \int_{\T} v(h(w))\, d\lambda(w) + 2\epsilon/3.
\end{align}
This implies that $f$ satisfies (a)-(e).

Next we use Lemma \ref{bigstep} for $f$, $E$ to get
analytic disc $g\in\cA_X$.
Property (i) implies that $g(z)\in W$ for every $z\in \T$. Since $\phi$
is continuous we may choose $\eta>0$ so small that property (ii) implies
 that
\begin{eqnarray}
\label{small-increase4}
\int_{\T\setminus E}\phi\circ g\,d\lambda\le \int_{\T\setminus E}\phi\circ f\,d\lambda+\epsilon/3,
\end{eqnarray}
and by (iii) we get  $g(0)=f(0)$.

Combining the inequalities (\ref{small-increase3}), (\ref{small-increase4}) and (d)
we obtain (\ref{eq:submean}).
This completes the proof of Theorem \ref{main}.\hfill \qed



\begin{thebibliography}{[AW]}

\bibitem[AG]{AG}
Andreotti, A.;  Grauert, H.:  
Th\'eor\`eme de finitude pour la cohomologie des espaces complexes.  
{Bull.\ Soc.\ Math.\ France}, \textbf{90}  193--259 (1962)

\bibitem[BS]{Bu-Schachermayer}
Bu,  S.\ Q.; Schachermayer, W.:
{Approximation of Jensen measures by image measures 
under holomorphic functions and applications},
Trans.\ Amer.\ Math.\ Soc., \textbf{331} 585--608 (1992)



\bibitem[Dr]{BD-Fourier}
Drinovec Drnov\v sek, B.:
On proper discs in complex manifolds.
{Ann.\ Inst.\ Fourier (Grenoble)}, \textbf{57} 1521--1535 (2007) 

\bibitem[DF1]{DF-Duke}
Drinovec Drnov\v sek, B.; Forstneri\v c, F.:
Holomorphic curves in complex spaces.
Duke Math.\ J., \textbf{139}, 203--254 (2007)


\bibitem[DF2]{DF-Indiana}
Drinovec Drnov\v sek, B.; Forstneri\v c, F.:
The Poletsky-Rosay theorem on singular complex spaces.
Indiana Univ.\ Math.\ J, in press. 
\texttt{http://arxiv.org/abs/1104.3968}

\bibitem[DF3]{DF-Polonici}
Drinovec Drnov\v sek, B.; Forstneri\v c, F.:
Disc functionals and Siciak-Zaharyuta extremal functions             %
on singular varieties.
Ann.\ Polon.\ Math., \textbf{106}, 171--191  (2012)




\bibitem[Fo1]{F-book}
Forstneri\v c, F.:
Stein Manifolds and Holomorphic Mappings
(The Homotopy Principle in Complex Analysis). 
Ergebnisse der Mathematik und ihrer Grenzgebiete, 
3.\ Folge, vol.\ 56, Springer-Verlag, Berlin-Heidelberg (2011)

\bibitem[Fo2]{FFAsian}
Forstneri\v c, F.:
Manifolds of holomorphic mappings from strongly pseudoconvex domains. 
Asian J.\ Math., \textbf{11},  113--126 (2007) 

\bibitem[Kl]{Klimek}
Klimek, M.: Pluripotential theory.
London Math.\ Soc.\ Monographs, new series, vol.\ 6,
The Clarendon Press, Oxford University Press, New York (1991)

\bibitem[Gr]{Grauert2} 
Grauert, H.:
Theory of $q$-convexity and $q$-concavity.
{\em Several complex variables, VII}, 259--284, 
\textit{Encyclopaedia Math.\ Sci.}, 74, 
Springer, Berlin, 1994.

\bibitem[Ku]{Kuzman}
Kuzman, U.:
Poletsky theory of discs in almost complex manifolds.
Comp.\ Var.\ Ell.\ Equat., in press. 

\bibitem[LP]{Larusson-Poletsky2012}
L\'arusson, F.; Poletsky, E.\ A.:
Plurisubharmonic subextensions as envelopes of disc functionals.
 Michigan Math.\ J., in press. 
\texttt{http://arxiv.org/abs/1201.5875}


\bibitem[LS1]{Larusson-Sigurdsson1998}
L\'arusson F.; Sigurdsson R.:
{Plurisubharmonic functions and analytic discs on manifolds},
J.\ Reine Angew.\ Math., \textbf{501},  1--39 (1998)

\bibitem[LS2]{Larusson-Sigurdsson2003}
L\'arusson F.; Sigurdsson R.:
{Plurisubharmonicity of envelopes of disc functionals on manifolds},
J.\ Reine Angew.\ Math., \textbf{555},  27--38 (2003)

\bibitem[LS3]{Larusson-Sigurdsson2005}
L\'arusson F.; Sigurdsson R.: 
{The Siciak-Zaharyuta extremal function as the envelope of disc functionals},
Ann.\ Polon.\ Math., \textbf{86},   177–-192 (2005)

\bibitem[LS4]{Larusson-Sigurdsson2007}
L\'arusson F.; Sigurdsson R.:
{Siciak-Zaharyuta extremal functions and polynomial hulls},
Ann.\ Polon.\ Math., \textbf{91},  235-–239 (2007)

\bibitem[LS5]{Larusson-Sigurdsson2009}
L\'arusson F.; Sigurdsson R.: 
{Siciak-Zaharyuta extremal functions, analytic discs 
and polynomial hulls},
Math.\ Ann., \textbf{345}, 159--174 (2009)







\bibitem[Ma]{Magnusson2011}
Magn\'usson, B. S.:
Extremal $\omega$-plurisubharmonic functions as envelopes of disc functionals.
Ark.\ Mat. \textbf{49}, 383–-399 (2011)

\bibitem[MS]{Magnusson-Sigurdsson2007}
Magn\'usson, B. S.; Sigurdsson R.: 
Disc formulas for the weighted Siciak-Zahariuta extremal function.
Ann.\ Polon.\ Math. \textbf{91}, 241–-247  (2007)


\bibitem[Po1]{Poletsky1991}
Poletsky, E.\ A.: 
{Plurisubharmonic functions as solutions of variational problems},
In: Several complex variables and complex geometry, Part 1 
(Santa Cruz, CA, 1989), 163--171,
Proc.\ Sympos.\ Pure Math., {52}, Part 1, 
Amer.\ Math.\ Soc., Providence, RI, 1991.

\bibitem[Po2]{Poletsky1993}
Poletsky, E.\ A.: Holomorphic currents. 
Indiana Univ.\ Math.\ J., \textbf{42}, 85--144 (1993)

\bibitem[RS]{Rashkovskii-Sigurdsson2005}
Rashkovskii, A.; Sigurdsson, R.:
Green functions with singularities along complex spaces. 
Internat.\ J.\ Math. \textbf{16}, 333-–355 (2005)

\bibitem[Ro1]{Rosay1}
Rosay, J.-P.:
Poletsky theory of disks on holomorphic manifolds. 
Indiana Univ.\ Math.\ J., \textbf{52}, 157--169 (2003)

\bibitem[Ro2]{Rosay2}
Rosay, J.-P.:
Approximation of non-holomorphic maps, and Poletsky theory of discs.
J.\ Korean Math. Soc., \textbf{40}, 423--434  (2003)


\end{thebibliography}
\end{document}